\def\a{{\alpha}}
\def\b{{\beta}}

\def\G{{\Gamma}}
\def\k{{\kappa}}

\def\p{{\pi}}

\def\cH{{\cal H}}

\def\cT{{\cal T}}

\def\pf{{\hfill$\Box$}}
\def\pr{{\prime}}
\def\proof{{\noindent {\it Proof}.\ \ }}
\def\ra{{\rangle}}
\def\la{{\langle}}

\def\l({{\left(}}
\def\r){{\right)}}

\def\({{\Biggl(}}
\def\){{\Biggr)}}
\def\[{{\Biggl[}}
\def\]{{\Biggr]}}

\documentclass[12pt]{article}

\newtheorem{thm}{Theorem}

\newtheorem{cnj}[thm]{Conjecture}
\newtheorem{cor}[thm]{Corollary}

\newtheorem{lem}[thm]{Lemma}

\newtheorem{res}[thm]{Result}

\usepackage{graphicx,amssymb,amsmath,verbatim}

\begin{document}

%
%
\title{Near-Universal Cycles for Subsets Exist}

\author{
Dawn Curtis\thanks{dawn.curtis@asu.edu},
Taylor Hines\thanks{taylor.hines@asu.edu},
Glenn Hurlbert\thanks{hurlbert@asu.edu},
Tatiana Moyer\thanks{tatiana.moyer@asu.edu}\\
    Department of Mathematics and Statistics\\
    Arizona State University,
    Tempe, AZ 85287-1804
}
\maketitle


\vspace{0.5 in}

%
%
\begin{abstract}
Let $S$ be a cyclic $n$-ary sequence. We say that $S$ is a {\it universal
cycle} ($(n,k)$-{\it Ucycle}) for $k$-subsets of $[n]$ if every such
subset appears exactly once contiguously in $S$, and is a {\it Ucycle
packing} if every such subset appears at most once. Few examples of
Ucycles are known to exist, so the relaxation to packings merits
investigation. A family $\{S_n\}$ of $(n,k)$-Ucycle packings for fixed $k$
is a {\it near-Ucycle} if the length of $S_n$ is $(1-o(1))\binom{n}{k}$.
In this paper we prove that near-$(n,k)$-Ucycles exist for all $k$.
\end{abstract}

\newpage

%
%
\section{Introduction}\label{Intro}

A {\it universal cycle (Ucycle) for $k$-subsets of $[n]$}, denoted $(n,k)$-UCS,
is a cyclic sequence of 
integers from $[n]=\{1,2,\ldots,n\}$ with the property that each $k$-subset
of $[n]$ appears exactly once consecutively in the sequence.
For example, $1234524135$ is a universal cycle for pairs of $[5]$.
Chung, et al. \cite{CDG}, defined universal cycles for a general class of 
combinatorial structures, generalizing both deBruijn sequences and Gray codes
(see \cite{Sav}).  The books \cite{Knu,Rus} contain a wealth of information
about generating such things efficiently.

A necessary condition for the existence of an $(n,k)$-UCS
is that $n\mid\binom{n}{k}$.
This is because symmetry demands that each symbol appears equally often.
The authors of \cite{CDG} conjectured that the necessary condition is
sufficient for large enough $n$ in terms of $k$ (some evidence \cite{Jac2} 
suggests $n\ge k+3$) and offered \$100 for its proof.  
\footnote{At the 2004
Banff Workshop on Generalizations of de Bruijn Cycles and Gray Codes it
was suggested by the second author that a modest inflationary rate should
revalue the prize near \$250.04.}

\begin{cnj}\label{CDG}\cite{CDG}
For all $k\ge 2$ there exists $n_0(k)$ such that for $n\ge n_0(k)$
Ucycles for $k$-subsets of $[n]$ exist if and only if $k$ divides 
$\binom{n}{k}$.
\end{cnj}

Progress on the conjecture has been slow. 
The $k=2$ case is trivial, corresponding to the existence of eulerian
circuits in $K_n$ if and only if $n$ is odd.
Jackson \cite{Jac1} proved the conjecture for $k=3$, and
constructed ucyles for $k=4$ and odd $n$, leaving 
the case $n\equiv 2\mod 8$ unresolved.
In \cite{Hur} we find the following result.

\begin{res}\label{OldRes}\cite{Hur}
Let $n_0(3)=8$, $n_0(4)=9$, and $n_0(6)=17$.  Then $(n,k)$-UCs exist for 
$k=$ 3, 4, and 6 with $n \geq n_0(k)$ and gcd$(n,k)=1$.
\end{res}

Note that $n_0(k)=3k$ suffices for $k\in \{3,4,6\}$.
It would be nice to lower $3k$ as much as possible.
To this end, Stevens et al. \cite{Ste} proved the following.

\begin{res}\label{stevens}\cite{Ste}
No $(k+2,k)$-UC exists for $k\ge 2$.
\end{res}

Combined with particular computer examples found by Jackson \cite{Jac2}
(e.g. $(n,k)=(10,4)$), this suggests that $n_0(k)=k+3$ may suffice.

In cases where no Ucycles have been found, or none exist, it is natural
to look for cycles with as many distinct subsets as possible; that is,
a {\it Ucycle packing}, such as the sequence
$$ S = 1345682 \  4678135 \  7812468 \  2345713 \  5678246 \  8123571 \  3456824 \  6781357 $$
for $(n,k)=(8,4)$.  Note that no $(8,4)$-UCS exists, and that $S$
accounts for 56 of the possible 70 subsets.
One might notice that these blocks shift upward by $3 \mod{8}$ from one
to the next.  This is important in the techniques that follow.

In their paper, Stevens, et al., show that the longest possible packing of a $(k+2,k)$-UCS
has length $k+2$ and a packing achieving this bound always exists.
Compared to the potential $\binom{k+2}{k}$ length, this shows that for 
$n=k+2$, we cannot get close to to a full universal cycle.  
To establish this notion formally, we define a {\it near-Ucycle} packing as a
sequence of Ucycle packings, one for each $n$, such that as $n \rightarrow \infty$,
asymptotically few $k$-subsets are omitted from the $(n,k)$-UCS packing. That is, 
the length of cycle $S_n$ is $(1-o(1))\binom{n}{k}$.
For example, if $n$ is even and $M$ is any perfect matching in $K_n$
then $K_n-M$ is eulerian.  In particular, any eulerian circuit is a
near-(n,2)-UCS of length $(1-\frac{1}{n-1})\binom{n}{2}$.
The purpose of this paper is to prove that near-Ucycles exist 
for all $k$.

\begin{thm}\label{Near}
For all $k$, near-$(n,k)$-Ucycles exist.
\end{thm}

We proceed in the proof of this theorem by analyzing the construction
of universal cycles.  As we will show, we can create a Ucycle packing
by selecting only those subsets that avoid certain structure.
We prove that there are asymptotically few such structured subsets.

%
%
\section{General Technique}\label{GenT}

\begin{figure}
\centerline{\includegraphics[height=1.75in]{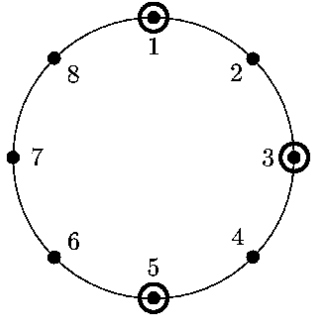}}
\caption{Form visualization of $\{1,3,5\}$}\label{135form}
\end{figure}

The general technique used to construct Ucycles originated from \cite{Jac1}.  It 
consists of classifying the component subsets by their structure and ordering 
them accordingly.  We write the $k$-subset $S$ of $[n]$ as $S = \{s_1, \dots, s_k\}$,
with $s_i < s_{i+1}$, and define the {\it form} of $S$ as $F = (f_1, \dots, f_k)$ by 
$f_i = s_{i+1}-s_i$, where indices are modulo $k$ and arithmetic is modulo $n$. 
That is, the form of a set is the ordered collection of distances between set elements. 
(see Figure \ref{135form})

Consider the following example.  For $(n,k)=(8,3)$, the set $\{1,3,5\}$ has form is 
$(2,2,4)$.  We consider the cyclic permutations of a form to be equivalent, so 
$(2,4,2)$ and $(4,2,2)$ can both serve as the form of $\{1,3,5\}$.  The choice 
of form has to do with how the form appears in the cycle.  For example, the form 
$F=(4,2,2)$ makes the sets $\{1,3,5\},\{2,4,6\},\dots,$ and $\{8,4,2\}$ appear as
$513, 624,\dots,$ and $482$, respectively.  Note that the last 2 in $F$ is not used for 
these sets, and so we may represent $F$ as $(4,2;2)$, or more simply $(4,2)$.
In our techniques below, we will restrict our attention to forms $(f_1,\dots,f_{k-1}; f_k)$
 having unique $f_k$.  In fact, we choose $f_k$ to be the largest unique entry. 

It is important to note that the method used here to construct the forms 
is somewhat limited.  As part of our definition, every form has entries 
whose sum is $n$.  We call such forms {\it simple}.  However, it is possible 
to relax this condition. For example, the set $\{5,3,1\}$ would have the 
form $(6,6,4)$.  If we allow the sum of form entries to be a multiple 
of $n$, we would have more freedom in representing the subsets, and 
hence a better method of constructing Ucycles could be developed.   
These forms we define as {\it crossing}.  In this paper, we will only use
simple forms.

The purpose of the forms is to model what occurs in a Ucycle, namely, that if 
$s_0s_1 \dots s_k$ appears in a Ucycle then the forms $(f_1,\dots,f_{k-1})$ and 
$(f_2,\dots,f_k)$ of the sets $\{s_0,\dots,s_{k-1}\}$ and $\{s_1,\dots,s_k\}$ overlap 
on $(f_2,\dots,f_{k-1})$.  This motivates the following definitions.

For a given form rep. $(f_1, \dots, f_{k-1})$, we define $(f_1, \dots, f_{k-2})$ to be its {\it prefix}, and $(f_{2}, \dots, f_{k-1})$ to be its {\it suffix}.

The {\it transition graph}, denoted $\cT_{n,k}$, is a graph whose vertices are the prefixes (and suffixes) of the form representations of $(n,k)$-UCS.  The directed edges are the form representations, drawn from prefix to suffix.  In order for this construction to produce a Ucycle, it is necessary that this transition graph be Eulerian.

Consider the cases of $(n,k)=(8,3)$ and $(n,k)=(10,4)$, and the forms of each:

We use these forms to construct $\mathcal{T}_{8,3}$ and $\mathcal{T}_{10,4}$, as shown in
Figure 2.

\begin{figure}
$$\begin{array}{ccc}
\includegraphics[height=1.5in]{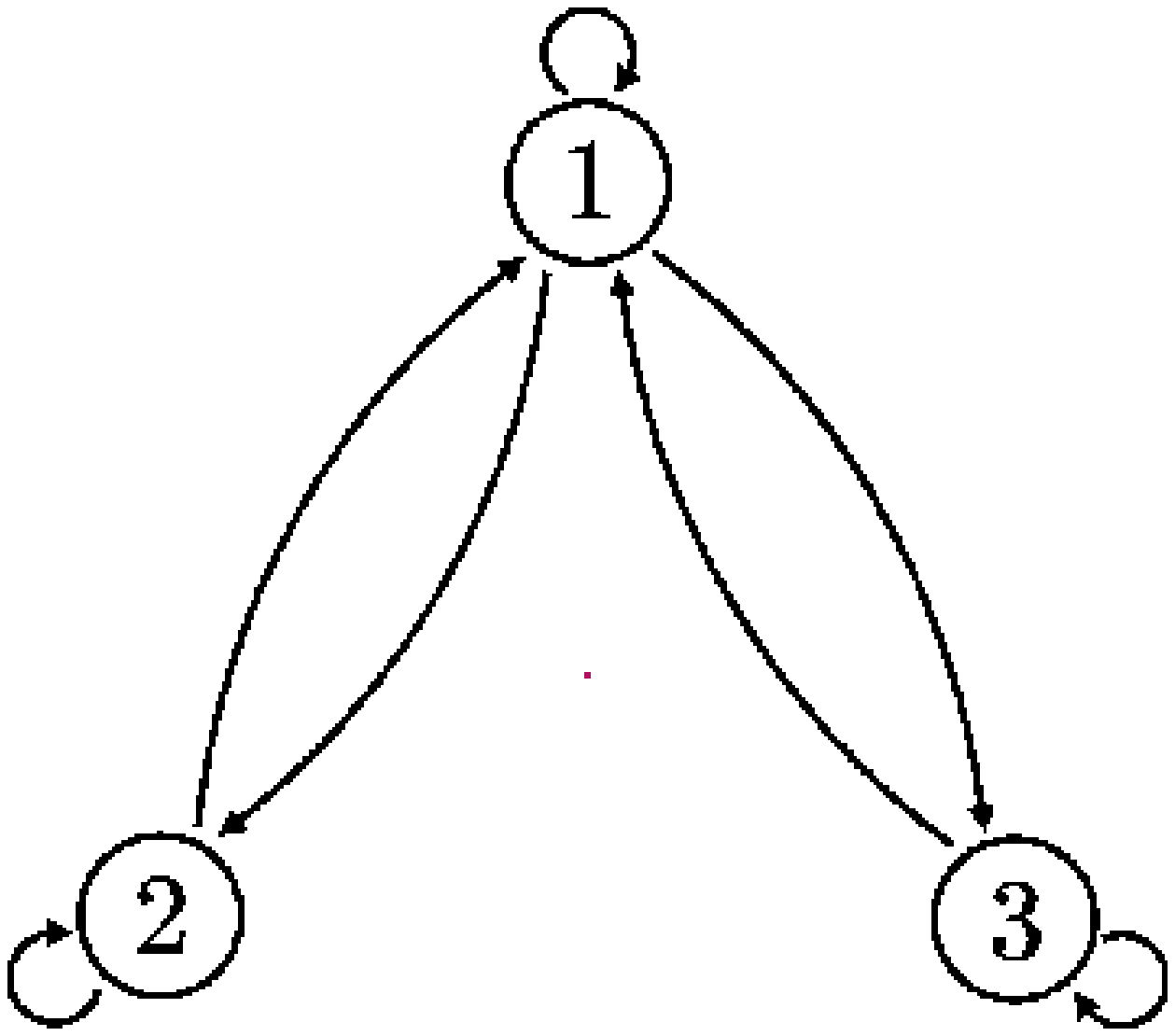}\label{T_83} &
\qquad&
\includegraphics[height=2in]{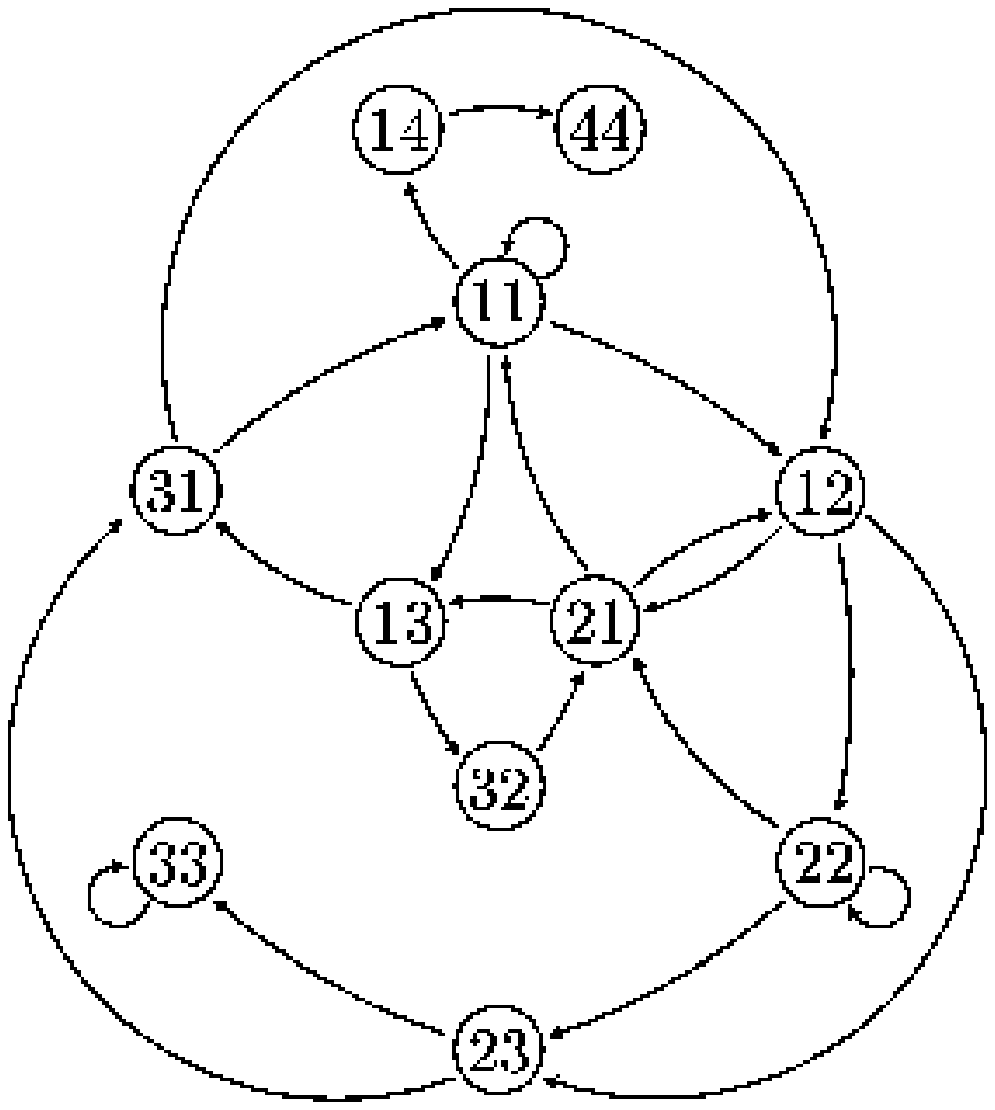}\label{T_104_bad} \\
&&\\
\cT_{8,3}&&\cT_{10,4}\\
\end{array}$$
\caption{Transition Graphs}
\end{figure}

The condition of evenness is highlighted by these two examples.  We know that a Ucycle is possible for $(n,k)=(8,3)$ but not for $(n,k)=(10,4)$.  As shown in the following result, this is directly connected to the fact that $\cT_{8,3}$ is Eulerian and $\cT_{10,4}$ is not.

\begin{res}\label{OldLem1}\cite{Hur}
If $\cT_{n,k}$ is Eulerian for some choice of form representations, then an $(n,k)$-UCS exists.
\end{res}

As we can see in $\cT_{10,4}$, the vertex $44$ has no out degree.  This is not necessarily the case.  
If we write the form $(1,1,4;4)$ as $(4,4,1;1)$, then we would instead connect $44 \rightarrow 41$.
Since this form has no unique entry, this is left ambiguous.  We call such forms \emph{bad}.
We define all forms with at least one unique element, a clear representative, as {\it good}.

Consider the previous example.  We can ignore all bad forms of $(n,k)=(10,4)$, and construct the transition graph, as shown in Figure 2.

\begin{figure}
\centerline{\includegraphics[height=2in]{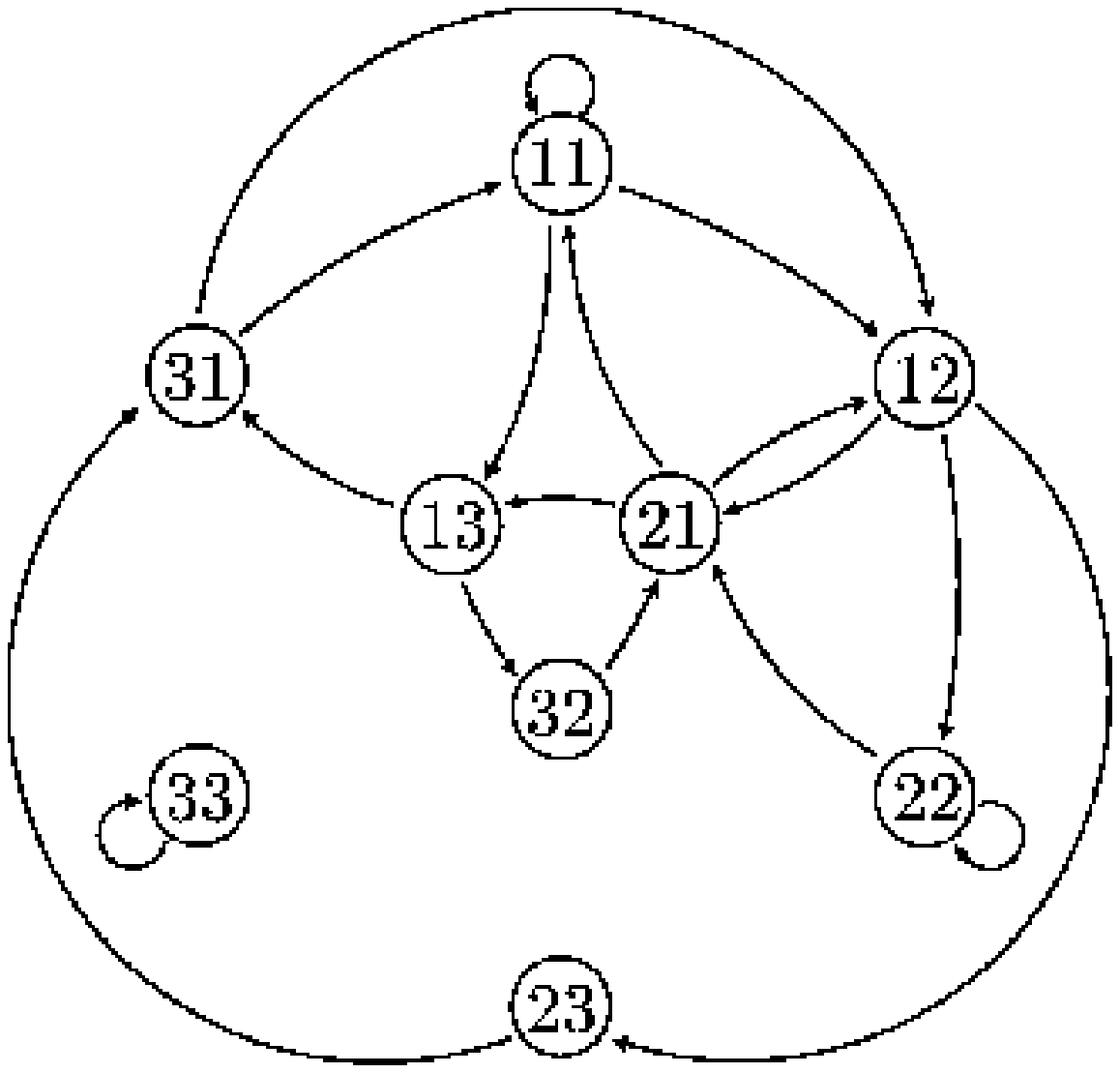}}\
\caption{Good $\cT_{10,4}$}
\label{T104}
\end{figure}

In this case, ignoring the bad forms yields an even transition graph.  
As stated in the following result, this is true for any connected transition 
graph.  Note, however, that it is not Eulerian since $((33))$ is not connected.
This problem disappears if we also ignore all such isolated cycles, 
as stated in the following result.  As we will show, the number of such isolated cycles
is negligible as $n \rightarrow \infty$.

\begin{res}\label{OldLem2}\cite{Hur}
If $\cH_{n,k}$ is connected and there are no bad forms for $k$-subsets of $[n]$, then $\cT_{n,k}$ is Eulerian.
\end{res}

We can construct a Ucycle packing by ignoring such bad forms for any $(n,k)$ pair
and restricting our attention to the largest component of $\cT_{n,k}$, the component containing
$((1 \dots 1))$.

%
%
\section{Main Result}\label{Main}

In order to prove that a near Ucycle packing is possible for any $k$, we must show that
an asymptotically large proportion of sets can be included in a Ucycle.
As we will show, the {\it good sets}, those belonging to good forms, 
can easily be included in a Ucycle packing, while the remaining {\it bad sets} cannot.
As we have indicated, we will only include good sets.  Of course, not all good sets
can be included.  It therefore remains to show that the good sets which can be included
asymptotically outnumber all other sets, and that a Ucycle packing that includes each of these sets can be created.
%
%
%
%
\subsection{Counting Good Sets}\label{CountGS}

Consider the example of the $(8,4)$-UCS.
We know that each form must have four entries, and, since each form is simple,
 the sum of these entries must be eight.  The number of times each entry appears 
 is also important, since the bad forms will have no entry appearing only once.
Therefore, we do not need all of its entries in order to determine whether or not a form is bad.
We need only to know how many times each entry appears.  
For example, the form $(1,2,2;3)$ has the unordered {\it pattern} $\la1,2,1\ra = \la2,1,1\ra$, 
essentially the list of multiplicities.  
\footnote{We can also denote the form $(1,2,2,3)$ as $(3^1 2^2 1^1)$, a notation we will find useful later}
Since this pattern contains a 1 as an entry, it has a unique element and therefore all 
forms with this {\it good pattern} are good.  In this case, the pattern entries define a 
partition of $4$, since each form has $4$ entries. In general, every pattern of an 
$(n,k)$-UCS is a partition of $k$.

Consider a {\it bad pattern} $P = \la p_1, \dots, p_t \ra$.  Clearly, $p_i \geq 2$ for 
all $1 \leq i \leq t$, so there exists a corresponding good pattern $P^{\pr} = \la p_1,\dots,p_t-1,1\ra$.
We define the function $\phi: \G \rightarrow \b$ from the set of of good patterns to the set of bad patterns
as $\p(P) = P^{\pr}$. Since this map applies to any bad pattern, we can see that $|\b| \leq |\G|$.
\footnote{Note that the number of good partitions of $k$ is equal to the number of partitions of $k-1$, denoted $p(k-1)$.
Since $p(k)\sim e^{\pi\sqrt{2k/3}}/4\sqrt{3}k\sim p(k-1)$,
almost all patterns are good.}

As we can see, many forms may belong to the same pattern.  In fact, the forms belonging 
to a particular pattern satisfy the equation $\sum_{i=1}^t p_ix_i = n$, for some positive integers
$x_1,\dots,x_t$, where $\la p_1,p_2,\dots,p_t\ra$ are the pattern entries.  
However, this equation imposes no order on the solution, and thus does not distinguish 
between two different forms that share the same entries.  We say that two such forms
belong to the same {\it class}.

For example, the forms $(1,2,2;3)$, $(2,1,2;3)$, and $(2,2,1;3)$ both share the class $[1,2,2;3]$. 
By convention, we order the entries of a class from smallest to largest entry.

As we know, the classes belonging to a pattern $P$ are all representations of $n$
as a positive distinct integer linear combination of $p_1, \dots, p_t$. 
We count the classes with the aid of Shur's theorem.
For $P=\la p_1,\dots,p_t\ra$, define the function
$$\psi(P)=\left| \left\{ (x_1,\dots,x_t) \mid \sum_{j=1}^{t} p_j x_j = n, 
x_j \geq 0 \right\} \right|\ .$$

\begin{thm}\label{Schur}\cite{Wilf}
Suppose that $gcd(p_1,\dots, p_t) = 1$, and define $k=\sum_{j=1}^{t} p_j$ and
$Q = \prod_{j=1}^{t} p_j$.  Then
$$\psi(P)\ \sim\ \frac{n^{t-1}}{(t-1)! Q}\ .$$
\end{thm}

Counting the classes requires a slightly more general result.  Namely, we want to extend Schur's theorem as follows.
We define
$$\psi^\prime(P)=\left| \left\{ (x_1,\dots,x_t) \mid\ \sum_{j=1}^{t} p_j x_j = n, 
\ x_j \geq 1, \ x_j \neq x_i \ \forall i \neq j \right\} \right|\ .$$

\begin{lem}\label{GenSchur}
Define $g=gcd(p_1,\dots, p_t)$, $k=\sum_{j=1}^{t} p_j$, and 
$Q =\prod_{j=1}^{t} p_j$.  Then
$$\psi^\prime(P)\ 
\sim\ \frac{(n-k)^{t-1}}{(t-1)! Qg^{t-1}} \ - 
\ \frac{(n-k)^{t-2}}{(t-2)!Qg^{t-2}} \sum_{i < j}\frac{p_i p_j}{p_i+p_j}\ .$$
\end{lem}

\proof

Since all form entries are positive integers, we count only integer solutions to $\sum_{j=1}^{t} p_j x_j = n$
such that $x_j \geq 1$.  This is equivalent to finding all integer solutions to the equation
$\sum_{j=1}^{t} p_j (x_j+1)$ $= \sum_{j=1}^{t} p_j x_j + k$ $= n$, without the positivity constraint.

We further modify the system to act for the condition $gcd(p_1,\dots, p_t) = g$.
The equation $\sum_{j=1}^{t} \frac{p_j}{g} x_j$ $= \sum_{j=1}^{t} q_j x_j$ $ = \frac{n-k}{g}$ 
is the same as our original equation, and clearly $gcd(q_1,\dots, q_t) = 1$.

Finally, we need to account for the fact that for a pattern of size $t$, each corresponding class
has exactly $t$ distinct entries.
If, for example, $P = \la 2,1,1\ra$, then $C = [1,1,3,5]$ is one possible class.  
However, $[3,3,3,1]$, although it is a valid solution to $\sum_{i=1}^{3}p_ix_i = n$, 
does not fulfill the requirement that $x_i \neq x_j$ for all $i \neq j$.
This is easily seen if we adopt the notation $C = [x_1^{p_1}, \dots, x_t^{p_t}]$.
Written in this way, $[1,1,3,5] = [1^2 3^1 5^1]$ and $[3,3,3,1] = [3^2 3^1 1^1]$
In the latter case, $x_1 = x_2$, so the class $[3^2 3^1 1^1]$ is not valid.
Therefore, we must account the number of such non-distinct solutions.  
If $x_i = x_k$, then the equation reduces to $(p_i+p_k)x_i + \sum_{j=1, j\neq i, j\neq k}^{t} p_j x_j =  n$.
We can then apply Schur's theorem to this equation as we did before.
For the new system, $Q^{\pr}$ $= (p_i+p_k)\prod_{j=1,j\neq i, j\neq k}^{t} p_j$ $= Q \frac{p_i + p_j}{p_i p_j}$, so the total number of non-distinct solutions to this equation is 
$\sim \frac{(n-k)^{t-2}}{(t-2)!Qg^{t-2}} \sum_{i < j}\frac{p_i p_j}{p_i+p_j}$
\pf

In order to count the classes, we now apply Lemma 8 to each pattern.
We want to show that as $n$ gets large, a given bad pattern $P$ has far fewer classes than its 
good component $\phi(P)$.
First, we must determine the size of $\sum_{i < j}\frac{p_i p_j}{p_i+p_j}$ for any pattern $P$.

For any $p_i$, we know that $\frac{p_i p_j}{p_i+p_j}$$= \frac{p_i(c-p_i)}{c}$ for some $c < k$,
which is maximized when $p_i = c/2$.  If $p_i = c/2$, then $\sum_{j=1}^t p_i$$= \frac{tc}{2}$$ = k $.
Thus $c = \frac{2k}{t}$.  Since each $p_i = c/2$, we maximize the value of this expression when $p_i = k/t$.  
Thus, $\sum_{i < j}\frac{p_i p_j}{p_i+p_j}$$\leq \binom{t}{2} \frac{(k/t)^2}{2k/t}$$\sim \frac{t^2}{2}\frac{k}{2t}$
$ = \frac{kt}{4} \leq k^2/4$.

Using this upper bound, we calculate that the number of classes belonging to a pattern $P$ is 

\begin{equation*}
c(P) \sim \displaystyle\frac{(n-k)^{t-1}}{(t-1)! Q g^{t-1}} - \displaystyle\frac{(n-k)^{t-2}}{(t-2)! Q g^{t-2}} \frac{k^2}{4} \sim \displaystyle\frac{(n-k)^{t-1}}{(t-1)! Q g^{t-1}} .
\end{equation*}

With this application of Schur's theorem, we can count the number of classes belonging to good patterns 
compared to the number of classes belonging to bad patterns as follows.
\begin{equation*}
\frac{c(P)}{c(\phi(P))} \sim 
\displaystyle\frac{(n-k)^{t-1}}{(t-1)! g^{t-1} \prod_{j=1}^{t} p_j}\displaystyle\frac{t! \prod_{j=1}^{t-1} p_j (p_t - 1)}{(n-k)^{t}}
\sim \frac{t(p_t-1)}{p_t n} \rightarrow 0
\end{equation*}
That is, the good classes asymptotically outnumber the bad classes.

It still remains to show that the number of forms of a bad class $C$ belonging to a bad pattern $P$
is no greater than the number of forms of a good class $C$ of $\phi(P)$.
As we know, the forms of a class are all permutations of the class entries modulo cyclic rotation.
Therefore, each class $C \in P$ has $\frac{(k-1)!}{\prod_{j=1}^{t}p_j!}$ forms, while
classes $C \in \phi(P)$ has $ \frac{k!}{\prod_{j=1}^{t-1}p_j! (p_t-1)}$ forms.
That is, the good forms outnumber the bad forms by a factor of $kp_t$.

Finally, it remains to count the sets.
By definition, all good forms have at least one unique element.  Thus the good sets have 
at least one unique difference between two elements.  This implies that every good form contains exactly $n$ sets.
\footnote{The bad forms do not have a unique element, so this is not always the case for bad sets.  
Instead, we only know that the number of sets contained in a bad form is a factor of $n$.
(Symmetry can reduce the number of sets; e.g. $\la 3 \ra$ has forms $(1,4,7)$, $(2,5,8)$, and $(3,6,9)$ when $n=9$.)}
Therefore, the number of sets per good form is at least the number of sets per bad form.

By counting the patterns, classes, forms and sets, one can see that almost every set is good.  
However, it remains to be shown that asymptotically many of these good sets can be arranged into a Ucycle packing.  
We will use the following lemmas to prove that this is true, and therefore a near Ucycle packing is always possible.

\begin{lem}\label{even}
If $\cT_{n,k}$ is restricted to the good classes then it is a union of cycles, and hence even.
\end{lem}

\proof
For a given class $C$ and transition graph $\cT_{n,k}$, we define the graph $\cT_{n,k}(C)$ to be the restriction of $\cT_{n,k}$ to the edges belonging to the forms of $C$.

If the class $C$ is good, then $\cT_{n,k}(C)$ is a cycle or union of cycles.
Each form of $C$ has a unique representative $c_k$, thus for any permutation $F$ of 
$\{c_1, \dots, c_{k-1}\}$, all cyclic permutations of $F$ are also forms of $C$.
Since $\cT_{n,k}(C)$ is a union of cycles for each good $C$, it is clear that $\cT_{n,k}$ restricted to the
good classes will be a union of Eulerian subgraphs, and is therefore even.
\pf

\begin{cor}
If the restriction of $\cT_{n,k}$ to good classes is connected, then it is Eulerian.
\end{cor}

\proof
Since the restriction of $\cT_{n,k}$ to good classes is a union of cycles, it is easily seen that 
if $\cT_{n,k}$ is connected, it is surely Eulerian.
\pf

By the result in \cite{Hur}, if $\cT_{n,k}$ is Eulerian, then an $(n,k)$-UCS exists.  Therefore, 
if we can show that the restriction of $\cT_{n,k}$ to good classes is connected, then it follows 
that a Ucycle packing exists that includes all good sets.  However, Fig. 2 demonstrates that this
is not always the case.  Instead, we prove that an asymptotically large component is connected.
Since we have proven that $(1-o(1))\binom{n}{k}$ sets are good, a Ucycle packing that includes 
all good sets is a near-packing. 

%
%
\subsection{Finding a Large Component}\label{LargeComp}

In order to study the components of the restriction of $\cT_{n,k}$ to good classes, we define the
{\it class graph}, denoted $\cH_{n,k}$, as the undirected graph whose vertices are all classes of 
$(n,k)$-UCS . An edge is drawn between the class representations that differ by only one entry. 
For example, $[1,2,2;5]$ and $[1,1,2;6]$ are connected in $\cH_{10,4}$.

If $C_1$ and $C_2$ are connected in $\cH_{n,k}$, then this means that 
$\cT_{n,k}(C_1)$ and $\cT_{n,k}(C_2)$ will also share vertices.
For example, $[1,2,2;5]$ and $[1,1,2;6]$ are connected in  $\cH_{10,4}$,
and correspond to the cycles
$$ ((12)) \  \rightarrow \  ((22)) \  \rightarrow \  ((21)) \  \rightarrow \  ((12)) $$
and
$$ ((11)) \  \rightarrow \  ((12)) \  \rightarrow \  ((21)) \  \rightarrow \  ((11)) \  .$$

Just because $C_1$ and $C_2$ are connected in $\cT{n,k}$, this does not guarantee 
that the union of $\cT_{n,k}(C_1)$ and $\cT_{n,k}(C_2)$ will be connected.  It could happen, for example,
that each $C_i$ has two components that connect, resulting in two components for $C_1 \cup C_2$.
However, as proven in \cite{Hur}, if $\cH_{n,k}$ is connected, then the union 
over all classes produces a connected $\cT_{n,k}$.  We will clarify this with the map $\kappa$,
defined as follows.

Let $C = [c_1^{p_1}, \dots, c_{t-1}^{p_{t-1}}; c_t]$ be a class, where the entry $c_i$ appears $p_i$ times, and $c_t$
is the largest singleton.
To connect the class good classes, we define the map $\k:c(\Gamma)\rightarrow c(\Gamma)$ by
$[c_1^{p_1}, \dots, c_{t-1}^{p_{t-1}}; c_t] \rightarrow [1,c_1^{p_1},\dots,c_{t-1}^{p_{t-1}-1};c_t+c_{t-1}-1]$
We then apply $\k$ again, each time adding another 1.

For example, $\k$ connects the class $[2,2,2;4]$ to the class $[1,1,1;7]$ of $\cH_{10,4}$ by the path
$$[2,2,2;4] \  \rightarrow \  [1,2,2;5] \  \rightarrow \  [1,1,2;6] \  \rightarrow \  [1,1,1;7] \  .$$

In this way, we are able to connect the majority of the classes to the class $[1,\dots,1;k-t+1]$.
However, using $\k$ does not work for every good class.
For example, the class $[3,3,3;1]$ is surely good.  However, $\k$ maps
$[3,3,3;1]$ to $[1,3,3;3]$; that is, to itself.  In general, if $c_t$ is the largest
singleton of a class $C$, then $c_t+c_{t-1}-1$ will be the largest singleton of $\k(C)$
only if $c_t > 1$.

%
%
\subsection{Counting Awesome Sets}\label{CountAS}

In order to circumvent the problem introduced above, we restrict our attention to the 
{\it awesome classes}, the good classes whose largest singleton is greater than 1.
Using Schur's theorem, we can show that the number of classes that are not awesome is 
negligible, as shown below.

$$\left| \left\{ (x_1,\dots,x_t) \mid \sum_{j=1}^{t}p_jx_j=n, \ 1=x_1=p_1<p_2<\dots<p_t,  \ x_j>1 \right\} \right|$$
\begin{eqnarray*}
 &=& \left| \left\{ (x_2,\dots,x_t) \mid \sum_{j=2}^{t}p_jx_j=n-1, \ 1<p_2<\dots<p_t,  \ x_j>1 \right\} \right|\\
&=& \left| \left\{ (y_2,\dots,y_t) \mid \sum_{j=2}^{t}p_jy_j=n-2k-1, \ 1<p_2<\dots<p_t,  \ y_j\geq0 \right\} \right|\\
&\sim& \frac{(n-2k-1)^{t-2}}{(t-2)!g^{t-2}Q}\\
&\ll& n^{t-1}
\end{eqnarray*}
for $n > 2k$.

Since the number of non-awesome classes is negligible compared to the number of total classes, we can
disregard them and restrict our attention to the awesome classes, which still comprise
an asymptotically large proportion of all subsets.

%
%
\subsection{Proof of Theorem \ref{Near}}\label{Proof}

Since almost all classes are awesome, we know that each awesome class is
connected to $[1,1,\dots,1;n-k+1]$ in $\cH_{n,k}$.
Since the awesome classes of $\cH_{n,k}$ are connected, the restriction of $\cT_{n,k}$
to the awesome classes is also connected.
If the restriction of $\cT_{n,k}$ to awesome classes is connected, then by lemma \ref{even}, we know
that all {\it awesome sets}, the sets belonging to awesome classes, can be connected to form
a Ucycle packing.  Since, as shown above, the awesome sets represent an asymptotically
large proportion of the total sets, it follows that this Ucycle packing is a near Ucycle packing.
\pf

%
%
\section{Remarks}\label{Remarks}

Many of the techniques presented here can be extended to other forms of Ucycle approximations.
For example, it is possible to include any set in a Ucycle packing by simply inserting the set elements anywhere in the
cycle.  We could therefore construct a {\it Ucycle covering} by simply appending all non-awesome sets onto
a near Ucycle created using the techniques we have described.
However, this is very inefficient because it increases the length of the cycle by $k$ for each added set 
instead of the desired 1.
In order to find a more elegant method of constructing Ucycle coverings, more complicated analysis is required.
As stated earlier, each bad form may only produce a factor of $n$ sets.  Therefore, the method used to connect awesome
sets into the Ucycle will not work, since traversing an Eulerian transition graph $n$ times is impossible for many bad forms.

One possible method of connecting bad sets in a Ucycle is to consider multiple form classes.
Currently, for a form $F = (f_1,\dots,f_k)$, we require $\sum_{i=1}^{k}f_i = n$.  However, it could be useful to consider
$\sum_{i=1}^{k}f_i = \a n$ for $\a > 1$.  This would allow much more freedom in representing sets, and therefore
more ways of connecting sets.

Finally, due to our proof that near Ucycles exist, we believe that we deserve asymptotically much of the prize money,
or $\$[1-o(1)](250.04)$.  Since we do not know the speed of the $o(1)$ term, we have made a conservative estimate
of $\$249.99$.

%
%
\bibliographystyle{plain}
%

%
%
\end{document}